\documentclass[12pt]{article}

\usepackage{amsmath}
\usepackage{amsthm}
\usepackage{amssymb}
\usepackage[mathscr]{eucal}

\theoremstyle{plain} 
\newtheorem{thm}{Theorem}
\newtheorem{prop}[thm]{Proposition}

\newtheorem{lem}[thm]{Lemma}

\theoremstyle{definition}

\newtheorem*{ques}{Question}

\theoremstyle{definition}
\newtheorem{rem}[thm]{Remark}
\newtheorem{rems}[thm]{Remarks}

\newcommand{\C}{\mathbb{C}}
\newcommand{\Q}{\mathbb{Q}}
\newcommand{\Qbar}{\overline{\mathbb{Q}}}
\newcommand{\Ql}{\mathbb{Q} _{\ell}}

\newcommand{\Z}{\mathbb{Z}}
\newcommand{\Zl}{\mathbb{Z} _{\ell}}

\newcommand{\Fp}{\mathbb{F} _p}
\newcommand{\Fq}{\mathbb{F} _q}
\newcommand{\Fpbar}{\overline{\mathbb{F} _p}}
\newcommand{\Fqbar}{\overline{\mathbb{F} _q}}

\newcommand{\alg}{\mathrm{alg}}

\newcommand{\End}{\mathrm{End}}

\newcommand{\Spec}{\mathrm{Spec}}

\newcommand{\GL}{\mathrm{GL}}

\renewcommand{\Sp}{\mathrm{Sp}}
\newcommand{\GSp}{\mathrm{GSp}}

\newcommand{\Fr}{\mathrm{Fr}}
\newcommand{\Frob}{\mathrm{Frob}}

\newcommand{\Tl}{\mathrm{T}_{\ell}}

\newcommand{\abs}[1]{\lvert#1\rvert}

\begin{document}

\title{A Note on the Existence of Absolutely Simple Jacobians}

\author{Ching-Li Chai and Frans Oort
\\Centre Emile Borel - UMS 839 IHP (CNRS/UPMC) - Paris}

\date{\small{version 2.1, 5/03/99}}

\maketitle


The purpose of this note is to answer a question to us
by M.~Sa\"{\i}di: 

\begin{ques} Let $g\ge 2$ be a positive integer.
Does there exist a projective nonsingular curve of genus $g$
over a finite field whose Jacobian is absolutely simple?
If so, are there infinitely many such curves?
\end{ques}
\smallskip

There is similar but simpler question, with
``finite fields'' replaced by ``number fields''.  
The answer for both questions is ``Yes'', see 
Remark \ref{three} (v).
\smallskip

It is easy to show that there exist infinitely many curves over
$\C$ whose Jacobian is simple:
the locus in the moduli space $\mathscr M_g$ of curves of 
genus $g$, consisting of those curves whose Jacobians are not simple, 
is the union of a countable
family of algebraic subvarieties of $\mathscr M_g$ of
lower dimension. So the union is not equal to 
$\mathscr M_g$.
\smallskip

Of course this naive argument does
not work for the countable fields $\Qbar$ and $\Fpbar$.
\smallskip

In this note we answer this question affirmatively.  
We show that for any family of curves $C\to S$ 
over a base scheme $S$ of
finite type over $\Fp$ 
(or over $\mathbb Q$) such that the monodromy is ``maximal'',
there exist infinitely many closed points $s\in S$
such that the Jacobian of $C_s$ is absolutely simple.
The proof uses results on $\ell$-adic 
representations, especially about Frobenius tori,
due to Serre, Chi, Larsen and Pink.
The proof also shows that the assumption on the
monodromy can be weakened; this is discussed in
Remark \ref{three} (i).
Also, the existence
of closed points $s\in \mathscr M_g(\Fpbar)$ with $\text{Jac}(C_s)$ absolutely
simple follows from a Chebotarev density argument, so 
the set of all such points $s$ has positive Dirichlet density
in the set $\abs{\mathscr M_g}$ of all closed points of 
$\mathscr M_{g,\Fp}$.
\smallskip

Here is a sketch of the basic argument. We are given
a family of abelian varieties $A\to S$ in characteristic $p$
with an assumption on its monodromy.
The goal is to show that there are many closed fibres in the
family which are absolutely simple.  
If for a closed point $s\in S$, the $\Ql$-Tate module 
$\Tl(A_{\bar s})\otimes_{\Zl}\Ql$  of the fibre
$A_s$ is $\Ql$-irreducible under the action of 
$\Fr_s^n$ for each power of the
Frobenius $\Fr_s$, then $A_s$ is absolutely simple.
Let $G^{\natural}$ be the set of all semisimple conjugacy classes of the
$\ell$-adic monodromy group $G$ for $A\to A$.
Let $Ir^{\natural}\subset G^{\natural}$  be the set of semisimple conjugacy
classes $[\gamma]$ such that $\Tl(A_{\bar s})\otimes_{\Zl}\Ql$ is 
irreducible as a $\Ql[\gamma]$ module.
If $Ir^{\natural}$ contains an open subset of $G^{\natural}$,
then we can apply the
Chebotarev density theorem to produce fibres $A_s$ in
the family which are absolutely simple.
In general, whether the subset $Ir^{\natural}$ is empty is a property of
the representation of $G$ on the Tate module.
In the exposition below, we assume for simplicity 
that the monodromy representation is 
the standard representation of the group of symplectic similitudes.
In this case $Ir^{\natural}$ contains the conjugacy classes of all 
generators of elliptic maximal tori, and Krasner's lemma
implies the openness we need.
\smallskip

\noindent{\bf Notations:} 
Let $S$ be a geometrically connected 
scheme smooth over a finite field $\Fq \supseteq \Fp$. 
Let $(\pi:A\to S, \lambda:A\to A^t)$
be a polarized abelian scheme over $S$ of relative 
dimension $g$. Let $\eta$ (resp.\ $\bar\eta$) 
be the generic point of $S$ (resp.\ a geometric generic point of $S$).
Let $\ell$ be a prime number different from $p$.
Let $\langle ,\rangle _{\lambda}$ be the symplectic pairing
on the Tate module $\Tl(A_{\bar\eta})\otimes \Ql$
attached to the polarization $\lambda$.
Let $*_{A,\lambda}$ be the symplectic
involution on the endomorphism algebra 
$\End(\Tl(A_{\bar\eta})\otimes \Ql)$ defined by 
$\langle ,\rangle _{\lambda}$.
Let $\Sp_{A, \lambda}$ (resp.\ $\GSp_{A, \lambda}$)
be the symplectic group (resp.\ the group of symplectic
similitudes) in $2g$ variables defined by $*_{A,\lambda}$.
Let $\rho:\pi_1(S, \bar\eta)\to 
\GSp_{A, \lambda}(\Q_{\ell})$ be the $\ell$-adic representations
attached to $\Tl(A_{\bar\eta})$.
Recall that the fundamental group $\pi_1(S,\bar\eta)$ fits 
canonically into an exact sequence
$$1
\longrightarrow \pi_1({\overline S},\bar\eta)
\longrightarrow \pi_1(S,\bar\eta)
\longrightarrow \text{Gal}(\Fqbar / \Fq)
\longrightarrow 1
\ \ ,
$$
where 
${\overline S} = S\times_{\Spec\,\Fq} \Spec\,\Fqbar$.
\smallskip

Let $G=G_A=\rho(\pi_1(S))$ 
(resp.\ $G_1=G_{A,1}=\rho(\pi_1({\overline S}, \bar\eta))\,$)
be the image of the fundamental group
(resp.\ the geometric fundamental group)
in $\GSp_{A,\lambda}(\Ql)$.  
Let $G^{\alg}=G_A^{\alg}$ 
(resp.\ $G_1^{\alg}=G_{A,1}^{\alg}$) be the
Zariski closure of $G_A$ (resp.\ $G_{A,1}$) in
$\GSp_{A, \lambda}$.
\smallskip

\noindent{\bf Construction/Definition:}
For each closed point $s\in S$, there is a Frobenius
element $\Frob_s$ in $\pi_1(S,\bar\eta)$
attached to $s$, well-defined up to conjugation.
The neutral component of the $\Ql$-Zariski closure of the
subgroup of $\GSp_{A, \lambda}$ generated by a
Frobenius element $\rho(\Frob_s)$ is a subtorus of
$\GSp_{A, \lambda}$; therefore it is also 
a subtorus of the Zariski closure $G_A^{\text{alg}}$.
This subtorus $T_s$ of $G_A^{\text{alg}}$ 
is well-defined up to conjugation in $G_A^{\text{alg}}(\Ql)$.
We call it ``the'' {\it Frobenius torus\/} attached to 
the closed point $s$, and denote it by $T_s$.
\smallskip

\begin{rem} The notion of Frobenius torus is due to Serre,
see \cite[\S3]{chi:thesis}, \cite[\S4, \S7]{larsen-pink:1992}.
The Frobenius torus $T_s$ can be defined over $\Q$; here we
use its ``$\ell$-adic realization''.
\end{rem}

In the next lemma we show that under suitable conditions 
on the monodromy of $A\to S$, there are many Frobenius tori $T_s$ which
are elliptic maximal tori; in other words $T_s$ is anisotropic
modulo the center of $\GSp_{2g,\lambda}$.
We denote by $\abs{S}$ the set of all closed points of $S$, and by
$\abs{S}_{\text{emFt}}$ the subset of all $s\in \abs{S}$ such that 
the Frobenius torus $T_s$ is an elliptic maximal torus 
of $\GSp_{A,\lambda}$ over $\Ql$.

\begin{lem} \label{ellip}
Notations as above.  Assume that the 
Zariski closure of the image
$\rho(\pi_1({\overline S},\bar\eta))$ of the
geometric fundamental group is equal to $\Sp_{A, \lambda}$.
Then $\abs{S}_{\text{emFt}}$ has positive Dirichlet
density in $\abs{S}$.  In particular $\abs{S}_{\text{emFt}}$
is an infinite set.
\end{lem}

\begin{proof}
This statement is certainly known.  We reproduce a
proof for the convenience of the reader.
As a first step, the assumption implies that $G$ is
an open subgroup in $G^{\alg}(\Ql)=\GSp_{A, \lambda}(\Ql)$.
Admit this for the moment.
\smallskip

The $\Ql$-rational maximal tori in $\GSp_{A, \lambda}$ are given by
maximal commutative semisimple subalgebras $B$ in
$\End(\Tl(A_{\bar\eta})\otimes\Ql)$ which are stable under
the involution $*_{A,\lambda}$.  For any such subalgebra $B$,
the $\Ql$-subspace of elements fixed under $*_{A,\lambda}$ has
dimension $g=\frac{1}{2} \dim_{\Ql}(B)$.
Among these maximal tori, the elliptic ones correspond to subfields of
dimension $2g$ over $\Ql$, and they are known to exist.
For any element $x\in \GSp_{A, \lambda}(\Ql)$, 
the neutral component of the Zariski closure of the
subgroup generated by $x$ is an elliptic maximal torus
if and only if the algebra 
$\Ql[x]\subset \End(\Tl(A_{\bar\eta})\otimes\Ql)$
is a field. 
Let $Em$ be the subset of $G$ consisting of all
elements of $x \in G$ such that 
$\Ql[x]\subset \End(\Tl(A_{\bar\eta})\otimes\Ql)$
is a field of degree $2g$ over $\Ql$.
By Krasner's lemma \cite[Prop.\ 3, p.\ 43]{lang:numb_th}, 
$Em$ is an open subset of $G$.
More concretely, if $\Ql[x_0]$ is a field of degree $2g$ over $\Ql$, then
for all $x\in \GSp_{A, \lambda}(\Ql)$ sufficiently close to $x_0$
in the $\ell$-adic topology, 
the characteristic polynomial and the eigenvalues of $x$ are close to 
those for $x_0$. 
So $\Ql[x]$ is a field of degree at least $2g$,  therefore the degree is
exactly $2g$.
It is known that $Em$ is not empty; 
so there exists an open normal subgroup $N\subset G$ such that
$Em$ is a union of cosets for $N$.
\smallskip

Recall that Chebotarev's density theorem states that
for any finite \'etale Galois covering $X\to Y$ between
normal integral schemes of finite type over $\Z$,
the Frobenius conjugacy classes in the covering group 
attached to the closed points of $Y$ is equidistributed,
see \cite[Thm.\ 7, p.\ 91]{serre:zetaL}.
Applying this density theorem to the finite \'etale
covering of $S$ corresponding 
to the finite quotient $G\to G/N$ of the monodromy group $G$, 
we see that there exist infinitely many Frobenius tori 
which are elliptic maximal tori in
$\GSp_{A,\lambda}$ over $\Ql$. 
This finishes the proof of Lemma \ref{ellip},
except for the fact that $G$ is open in $G^{\alg}(\Ql)$.
This is done in the next two paragraphs. 
The reader may wish to skip them because it is a little
technical and not central to this note; cf.\ Remark \ref{two} (i).
\smallskip

Both $G$ and $G_1$ are closed subgroups of
$\GSp_{A,\lambda}(\Ql)$; hence they are analytic Lie
subgroups of $\GSp_{A,\lambda}(\Ql)$.
Let $\mathfrak g$
(resp.\ $\mathfrak g_1$)  be the Lie algebra of
$G$ (resp.\ $G_1$),  and let 
$\mathfrak g^{\text{alg}}$ 
(resp.\ $\mathfrak g_1^{\text{alg}}$) be the Lie algebra  
of $G^{\alg}$ (resp.\ $G_1^{\alg}$).
The assumption of the theorem means that 
$\mathfrak g_1^{\text{alg}}$ 
is equal to the Lie algebra $\mathfrak{sp}_{2g,\lambda}$ of
$\Sp_{A,\lambda}$.  
Therefore $\mathfrak g^{\text{alg}}$ 
is equal to the Lie algebra $\mathfrak{gsp}_{2g,\lambda}$ of
$\GSp_{A,\lambda}$.
\smallskip

By \cite[\S7, Cor.\ 7.9]{borel:alg_grp}, 
we have 
$
[\mathfrak g_1^{\text{alg}}, \mathfrak g_1^{\text{alg}}]
\subseteq \mathfrak g_1
$.
We reproduce a proof here. 
The algebraic subgroup of $\GSp_{A,\lambda}$ consisting of
all elements $x\in \GSp_{A,\lambda}$ such that
$\text{Ad}(x)(\mathfrak g_1)\subseteq \mathfrak g_1$
clearly contains $G$, hence 
$
[\mathfrak g_1^{\text{alg}}, \mathfrak g_1]
\subseteq \mathfrak g_1
$. 
Repeating this argument again, one concludes that 
$
[\mathfrak g_1^{\text{alg}}, \mathfrak g_1^{\text{alg}}]
\subseteq \mathfrak g_1$.
\smallskip

From $[\mathfrak g_1^{\text{alg}}, \mathfrak g_1^{\text{alg}}]
\subseteq \mathfrak g_1$, we deduce that 
$\mathfrak g_1=\mathfrak{sp}_{2g,\lambda}$ and
$\mathfrak g=\mathfrak{gsp}_{2g,\lambda}$.
Especially $G$ is an open subgroup of $\GSp_{A,\lambda}(\Ql)$.
\end{proof}

\begin{rems}\label{two} (i) The argument to show 
that $\mathfrak g_1$ is big
since it contains the derived algebra of its algebraic envelope
appeared in \cite{bogomolov:l-adic}.


(ii) It is easy to find examples which satisfy the assumption
in Lemma \ref{ellip}. For instance, one can take $S$ to be
equal to the moduli space $\mathscr M_{g,n}$ of curves of
genus $g$ over $\Fp$ with the principal symplectic level-$n$ structure, or a 
generic curve in $\mathscr M_g$, with $n\ge 3$ and $(n, p\ell)=1$.
In this case the monodromy group $G$ is equal to $\GSp_{2g}(\Ql)$;
this is equivalent to
\cite[Thm.\ 5.15, p.\ 108]{deligne-mumford}, 
that the moduli spaces $\mathscr M_{g, n\ell^i}$ 
are geometrically connected.
Then Proposition \ref{simp}
below answers the Question at the beginning of this note.
\end{rems}

\begin{prop}\label{simp} Notation and assumptions as in Lemma \ref{ellip}.
Then for each $s\in \abs{S}_{\text{emFt}}$, 
the fiber $A_s$ is an absolutely simple abelian variety.
\end{prop}

\begin{proof}
Suppose that $s\in S$ is a closed point such that the 
Frobenius torus $T_s$ is an elliptic maximal torus
in $\GSp_{A,\lambda}$ over $\Ql$.
We know that each elliptic maximal torus $T$ in
$\GSp_{A,\lambda}$ is given by a subfield $F$ in
$\End(\Tl(A_{\bar\eta})\otimes\Ql)$ of degree $2g$ stable under
$*_{A,\lambda}$, hence $T$ acts irreducibly on 
$\Tl(A_{\bar\eta})\otimes\Ql$.
So the action of any open subgroup of $\text{Gal}(\bar s/s)$ on
the $\ell$-adic Tate module of $A_s$ is irreducible over
$\Ql$.
This implies that $A_s$ is not isogenous to a non-trivial product
over any finite extension of the residue field 
$\kappa(s)$ of $s$.
\end{proof}

\begin{rems}\label{three} 
(i) In the proof above, what we really need about 
the monodromy of $\pi:A\to S$ is that the neutral component
$G^{\text{alg}}$ of the Zariski closure of $\rho(\pi_1(S))$
has an elliptic maximal torus $T$ over $\Ql$ such that
the restriction to $T$ of the representation of $G^{\text{alg}}$
on the $\ell$-adic Tate module is irreducible.
For instance, if $G^{\text{alg}}$ is isomorphic to the derived group
$(\mathrm{Res}_{E/\Ql}\GL_2)^{\text{der}}$
of $\mathrm{Res}_{E/\Ql}\GL_2)$ for
a finite extension $E$ of $\Ql$
and the representation 
is isomorphic to the standard representation
of $\mathrm{Res}_{E/\Ql}\GL_2)$, 
the argument still works and one gets absolutely
simple fibers in the abelian scheme in question.

(ii) On the other hand, 
suppose $\pi:A\to S$ is the universal family of 
abelian surfaces over the Shimura curve $S$ 
attached to a quaternion division algebra $B$
over $\Q$, such that $B$ is split at $p$, 
and $\ell$ is a prime number different from $p$.
Then the representation of every Frobenius torus $T_s$
on the $\ell$-adic Tate module is reducible.
The group $G^{\text{alg}}$ is the group
attached to $B^{\text{opp},\times}$ where $B^{\text{opp}}$ 
is the opposite algebra of $B$, 
and the representation on the $\ell$-adic Tate module
comes from the regular representation of $B^{\text{opp}}$.
In this example one can show that none of the fibres $A_s$
is absolutely simple.

(iii) If the representation of the Frobenius torus $T_s$ on
the $\ell$-adic Tate module is reducible, 
then $\End(A_{\bar s})\otimes \Ql$ contains a 
nontrivial idempotent.
But this in general does not imply that 
$\End(A_{\bar s})\otimes \Q$ contains a nontrivial idempotent, so
$A_s$ over the closed point $s$ may or may not be absolutely simple.

(iv) If $A_K$ is an abelian variety of dimension $g$ 
over a number field 
such that the Zariski closure of the image of the $\ell$-adic
Galois representation attached to $A$ is equal to $\GSp_{2g}$,
the same argument shows that there exist finite
places $v$ of $K$ such that the reduction $A_v$ of $A$ at $v$ is 
an absolutely simple abelian variety; 
such places $v$ form a subset of positive density.
One can also change $\GSp_{2g}$ to a group satisfying the
property specified in (i).

By \cite[1.4]{noot:ord_red}, due to Serre, 
for any abelian variety $A_L$ over a field $L$ of finite type
over $\Q$, there exist specializations $B_K$ to number fields
$K$, such that the images of the $\ell$-adic 
Galois representations for $A_L$ and $B_K$ are equal.
Therefore there exist examples satisfying the assumption in 
the previous paragraph.

(v) Proposition \ref{simp} gives an
affirmative answer to the Question at the
beginning of this note.  
The similar question, with $\Fpbar$ replaced by $\Qbar$,
was answered by \cite[Cor.\ 1.5]{noot:ord_red}.
and also by \cite[Thm.\ 0.6.1, p.\ 8]{andre:motive}.
One can also deduce it from the $\Fpbar$-version because of 
Lemma \ref{qbar} below. 

\end{rems}

\begin{lem}\label{qbar}
Suppose that $A_K$ is an abelian variety over a number
field $K$ with good reduction at a finite place $v$ of $K$.
If the reduction $A_v$ of $A$ is absolutely simple, then the
abelian variety $A$ is also absolutely simple.
\end{lem}

\begin{proof}
Since the hypothesis does not change if one extends the base 
field $K$ to a bigger number field $L$,
it suffices to prove that the abelian variety $A$ is simple over $K$.
Recall that whether the abelian variety $A$ is simple or not
is a question on the existence of an endomorphism
$T$ of $A$ such that $T^2=n\cdot T$ for some non-zero integer $n$,
while $\text{Ker}(T)$ is not finite.  
Since every endomorphism of $A$ over $K$ extends to an 
endomorphism of the abelian scheme $A$ over $\mathscr O_{K,v}$
which extends $A$, we are done.
\end{proof}

\noindent{\it Acknowledgement\/}: The authors would like to
thank B.~Moonen, R. Noot and Y.~Zarhin for suggestions.
The work was done at Institut Henri Poincar\'e
during the 1999 Diophantine Geometry Trimester.
We thank Centre Emile Borel for hospitality and for
providing a stimulating environment.


\vspace{10mm}
\begin{tabbing} Ching-Li Chai \\ Department of Mathematics
\\University of Pennsylvania \\ 209 South 33rd St.
\\Philadelphia, PA 19104-6395 \\ email: chai@math.upenn.edu
\end{tabbing}

\vspace{10mm}
\begin{tabbing} Frans Oort  \\ Mathematisch Instituut  \\ Budapestlaan 6  
\\ NL: 3508 TA Utrecht  \\ The Netherlands \\ email: oort@math.uu.nl  \\
\end{tabbing}

\end{document}